\pdfoutput=1
\documentclass{amsart}
\usepackage[utf8]{inputenc}

\usepackage[margin=1in]{geometry}
\usepackage[expansion=false]{microtype}

\usepackage{amsmath, amsfonts, amssymb, amsthm, amsopn}
\usepackage{eucal}

\usepackage{tikz}
\usetikzlibrary{diagrams}

\usepackage[pdfusetitle,unicode,hidelinks]{hyperref}

\newcommand{\ie}{{\sfcode`\.1000 i.e.}}

\numberwithin{equation}{section}

\theoremstyle{plain}
\newtheorem*{theorem*}{Theorem}

\newtheorem{theorem}[equation]{Theorem}
\newtheorem{proposition}[equation]{Proposition}
\newtheorem{lemma}[equation]{Lemma}
\newtheorem{corollary}[equation]{Corollary}

\theoremstyle{definition}
\newtheorem{definition}[equation]{Definition}

\theoremstyle{remark}
\newtheorem{remark}[equation]{Remark}

\let\scr=\mathcal
\let\bb=\mathbb

\def\Z{\bb Z}
\def\Q{\bb Q}

\def\1{\mathbf 1}
\def\T{\bb T}

\let\del=\partial

\let\into=\hookrightarrow
\let\onto=\twoheadrightarrow
\let\tens=\otimes

\def\ph{\mathord-}

\def\abs#1{\lvert #1\rvert}

\def\id{\mathrm{id}}
\DeclareMathOperator{\Hom}{Hom}

\DeclareMathOperator{\Fun}{Fun}

\DeclareMathOperator{\Aut}{Aut}
\DeclareMathOperator{\Map}{Map{}}
\DeclareMathOperator{\ob}{ob{}}
\DeclareMathOperator{\Spec}{Spec}

\DeclareMathOperator{\Tot}{Tot}

\def\PSh{\mathrm{PSh}{}}

\def\op{\mathrm{op}}

\def\Cat{\scr C\mathrm{at}}

\def\HH{\mathrm{HH}}
\def\HC{\mathrm{HC}}
\def\HN{\mathrm{HN}}
\def\HP{\mathrm{HP}}
\def\CC{\mathrm{CC}}
\def\CN{\mathrm{CN}}
\def\CP{\mathrm{CP}}
\def\BC{\mathrm{BC}}
\def\BN{\mathrm{BN}}
\def\BP{\mathrm{BP}}
\def\CAlg{\mathrm{CAlg}{}}

\def\Mod{\mathrm{Mod}{}}
\def\Comod{\mathrm{Comod}{}}
\def\Ch{\mathrm{Ch}{}}

\def\LL{\tilde\Lambda}

\let\lim=\relax
\DeclareMathOperator*{\lim}{lim}
\DeclareMathOperator*{\colim}{colim}

\title{The homotopy fixed points of the circle action on Hochschild homology}
\author{Marc Hoyois}
\date{\today}

\begin{document}

\begin{abstract}
	We show that Connes' $B$-operator on a cyclic differential graded $k$-module $M$ is a model for the canonical circle action on the geometric realization of $M$. This implies that the negative cyclic homology and the periodic cyclic homology of a differential graded category can be identified with the homotopy fixed points and the Tate fixed points of the circle action on its Hochschild complex.
\end{abstract}

\maketitle

Let $k$ be a commutative ring and let $A$ be a flat associative $k$-algebra. The Hochschild complex $\HH(A)$ of $A$ with coefficients in itself is defined as the normalization of a simplicial $k$-module $A^\natural$ with
\[
A^\natural_n=A^{\tens_k (n+1)}.
\]
The simplicial $k$-module $A^\natural$ is in fact a \emph{cyclic} $k$-module: it extends to a contravariant functor on Connes' cyclic category $\Lambda$. As we will see below, it follows that the chain complex $\HH(A)$ acquires a canonical action of the circle group $\T$.
The cyclic homology, negative cyclic homology, and periodic cyclic homology of $A$ over $k$ are classically defined by means of explicit bicomplexes.
The goal of this note is to show that:
\begin{enumerate}
	\item The \emph{cyclic homology} $\HC(A)$ coincides with the \emph{homotopy orbits} of the $\T$-action on $\HH(A)$.
	\item The \emph{negative cyclic homology} $\HN(A)$ coincides with the \emph{homotopy fixed points} of the $\T$-action on $\HH(A)$.
	\item The \emph{periodic cyclic homology} $\HP(A)$ coincides with the \emph{Tate fixed points} of the $\T$-action on $\HH(A)$.
\end{enumerate}
The first statement is due to Kassel \cite[Proposition A.5]{Kassel}. The second statement is well-known to experts, but a proof seems to be missing from the literature. This gap was mentioned in the introductions to \cite{TV} and \cite{TV2}, and it was partially filled in \cite{TV}, where (2) is proved at the level of connected components for $A$ a smooth commutative $k$-algebra and $\Q\subset k$. 

We will proceed as follows:
\begin{itemize}
	\item In \S\ref{sec:general}, we recall abstract definitions of cyclic, negative cyclic, and periodic cyclic homology in a more general context, namely for $\infty$-categories enriched in a symmetric monoidal $\infty$-category. 
	\item In \S\ref{sec:concrete}, we show that in the special case of differential graded categories over a commutative ring, the abstract definitions recover the classical ones.
\end{itemize}

\subsubsection*{Acknowledgments} I thank Thomas Nikolaus for spotting an error in the proof of Proposition \ref{prop:realization} in a previous version of this note.

\section{Cyclic homology of enriched \texorpdfstring{$\infty$}{∞}-categories}
\label{sec:general}

Let $\scr E$ be a presentably \emph{symmetric} monoidal $\infty$-category, for instance the $\infty$-category $\Mod_k$ for $k$ an $E_\infty$-ring. We denote by $\Cat(\scr E)$ the $\infty$-category of categorical algebras in $\scr E$ in the sense of \cite{GH}, \ie, $\scr E$-enriched $\infty$-categories with a specified space of objects $\ob(\scr C)$. 
Let $\Lambda$ be Connes' cyclic category \cite{Connes}, with objects $[n]$ for $n\in\mathbb N$. 
 We can associate to every $\scr C\in \Cat(\scr E)$ a cyclic object $\scr C^\natural\colon \Lambda^\op \to \scr E$ with
\[
\scr C^\natural_n = \colim_{a_0,\dotsc,a_n \in \ob(\scr C)} \scr C(a_0,a_1)\otimes \dotsb \otimes\scr C(a_{n-1},a_n) \otimes \scr C(a_n,a_0).
\] 
The cyclic category $\Lambda$ contains the simplex category $\Delta$ as a subcategory, and the \emph{Hochschild homology} of $\scr C$ (with coefficients in itself) is the geometric realization of the restriction of $\scr C^\natural$ to $\Delta^\op$:
\[
 \HH(\scr C)=\colim_{n\in\Delta^\op}\scr C^\natural_n\in\scr E.
\]
We refer to \cite{AMGR} for a precise construction of $\scr C^\natural$ in this context, and for a proof that $\HH(\scr C)$ depends only on the $\scr E$-enriched $\infty$-category presented by $\scr C$ (\ie, the functor $\HH$ inverts fully faithful essentially surjective morphisms).

We now recall the canonical circle action on $\HH(\scr C)$.
Let $\Lambda\to \LL$ denote the $\infty$-groupoid completion of $\Lambda$, and let $\T$ be the automorphism $\infty$-group of $[0]$ in $\LL$. Since $\Lambda$ is connected, there is a canonical equivalence $B\T\simeq\LL$.
 Let $\PSh(\Lambda,\scr E)$ be the $\infty$-category of $\scr E$-valued presheaves on $\Lambda$, and let $\PSh_\simeq(\Lambda,\scr E)\subset\PSh(\Lambda,\scr E)$ be the full subcategory of presheaves sending all morphisms of $\Lambda$ to equivalences. 
We have an obvious equivalence
\[
\PSh_\simeq(\Lambda,\scr E)\simeq \PSh(B\T,\scr E).
\]
Since $\scr E$ is presentable and $\Lambda$ is small, $\PSh_\simeq(\Lambda,\scr E)$ is a reflective subcategory of $\PSh(\Lambda,\scr E)$. We denote by
\[
\abs\ph\colon \PSh(\Lambda,\scr E)\to \PSh_\simeq(\Lambda,\scr E)\simeq \PSh(B\T,\scr E)
\]
the left adjoint to the inclusion.
The morphisms
\[
*\xrightarrow{i} B\T\xrightarrow{p} *
\]
each induce three functors between the categories of presheaves. We will write
\[u_\T=i^*,\quad (\ph)_{h\T}=p_!,\quad (\ph)^{h\T}=p_*\]
for the forgetful functor, the $\T$-orbit functor, and the $\T$-fixed point functor, respectively.

\begin{proposition}
	\label{prop:realization}
	Let $\scr E$ be a presentable $\infty$-category and let $X\in\PSh(\Lambda,\scr E)$ be a cyclic object. There is a natural equivalence \[u_\T\abs{X}\simeq\colim_{[n]\in\Delta^\op}X([n]).\]
\end{proposition}

\def\Aut{\mathrm{Aut}}

\begin{proof}
	Let $j\colon \Delta\to\Lambda$ be the inclusion.
	Let $X\in\PSh_\simeq(\Delta,\scr E)$ and let $j_*(X)\in\PSh(\Lambda,\scr E)$ be the right Kan extension of $X$.
	We claim that $j_*(X)$ inverts all morphisms in $\Lambda$. By the formula for right Kan extension, it suffices to show that following: for every $[n]$, the functor of comma categories $\Delta\times_{\Lambda}\Lambda_{/[n]} \to \Delta\times_{\Lambda}\Lambda_{/[0]}$ induced by the unique map $[n]\to[0]$ in $\Delta$ is a weak equivalence.
	Recall that every morphism in $\Lambda$ can be written uniquely as a composite $h\circ t$ where $t$ is an automorphism and $h$ is in $\Delta$ \cite[Theorem 6.1.3]{Loday}.
	If $\Gamma = \Delta\times_{\Lambda}\Lambda_{/[0]}$, an object of $\Gamma$ is a pair $([m], t)$ with $[m]\in\Delta$ and $t \in \Aut_\Lambda([m]) = C_{m+1}$, and a morphism $([m],t) \to ([m'], t')$ is a map $h\colon [m] \to [m']$ in $\Delta$ such that $t'ht^{-1}$ is in $\Delta$. It is then clear that the functor 
	\[
	\Gamma\to\Delta, \quad ([m],t)\mapsto [m],\quad h\mapsto t'ht^{-1},
	\]
	is a cartesian fibration (whose fibers are sets). Moreover, the functor $\Delta\times_{\Lambda}\Lambda_{/[n]} \to \Delta\times_{\Lambda}\Lambda_{/[0]}$ can be identified with the projection $\Gamma \times_{\Delta} \Delta_{/[n]} \to\Gamma$.
	Since $\Delta$ is cosifted, the forgetful functor $\Delta_{/[n]}\to\Delta$ is coinitial. The pullback of a coinitial functor along a cartesian fibration is still coinitial \cite[Proposition 4.1.2.15]{HTT}, so $\Gamma \times_{\Delta} \Delta_{/[n]} \to\Gamma$ is coinitial and in particular a weak equivalence, as desired.
	
	Thus, we have a commuting square
	\begin{tikzmath}
		\diagram{\PSh_\simeq(\Delta,\scr E) & \PSh(\Delta,\scr E) \\
		\PSh_\simeq (\Lambda,\scr E) & \PSh(\Lambda,\scr E)\rlap. \\};
		\arrows (11-) edge[c->] (-12) (21-) edge[c->] (-22)
		(11) edge node[left]{$j_*$} (21) (12) edge node[right]{$j_*$} (22);
	\end{tikzmath}
	Since $\Delta^\op$ is weakly contractible, evaluation at $[0]$ is an equivalence $\PSh_\simeq(\Delta,\scr E)\simeq \scr E$. The left adjoint square, followed by evaluation at $[0]$, says that $u_\T\abs{\ph}\simeq\colim j^*(\ph)$, as desired.
\end{proof}

It follows from Proposition~\ref{prop:realization} that $u_\T\abs{\scr C^\natural}\simeq\HH(\scr C)$, so that $\HH(\scr C)$ acquires a canonical action of the $\infty$-group $\T$. 
As another corollary, we recover the following computation of Connes \cite[Théorème 10]{Connes}:

\begin{corollary}\label{cor:KZ2}
	$\LL\simeq K(\Z,2)$. 
\end{corollary}

\begin{proof}
	If $X\in\PSh_\simeq(\Lambda)$, then, by Yoneda, $\Map(\Lambda^0,X)\simeq \Map(\LL^0,X)$. In other words, the canonical map $\Lambda^0\to\LL^0$ induces an equivalence $\abs{\Lambda^0}\simeq \LL^0$, and hence $u_\T\abs{\Lambda^0}\simeq\T$.
	On the other hand, the underlying simplicial set of $\Lambda^0$ is $\Delta^1/\del\Delta^1$, so $u_\T\abs{\Lambda^0}\simeq K(\Z,1)$ by Proposition~\ref{prop:realization}. Thus, $\T$ is a $K(\Z,1)$, and hence $B\T\simeq\LL$ is a $K(\Z,2)$.
\end{proof}

In particular, $\T$ is equivalent to the circle as an $\infty$-group, which justifies the notation.
If $\scr E$ is stable, Atiyah duality for the circle provides the \emph{norm map} $\nu_{\T}\colon \Sigma^{\mathfrak t}(\ph)_{h\T} \to (\ph)^{h\T}$, where $\mathfrak t$ is the Lie algebra of $\T$ and $\Sigma^{\mathfrak t}$ is suspension by its one-point compactification. Explicitly, if $E\in\scr E$ has an action of $\T$, the norm is induced by the $(\T\times\T)$-equivariant composition
\[
\Sigma^{\mathfrak t}E \to \Sigma^{\mathfrak t}\Hom(\Sigma^\infty_+\T, E) \simeq \Sigma^{\mathfrak t}(\Sigma^\infty_+\T)^\vee \otimes E \simeq \Sigma^\infty_+\T \otimes E\to E,
\]
where the first map is the diagonal, the third is Atiyah duality, and the last is the action. 
The cofiber of $\nu_\T$ is the \emph{Tate fixed point functor} $(\ph)^{t\T}$.

\begin{definition}
	Let $\scr E$ be a presentably symmetric monoidal $\infty$-category and let $\scr C\in\Cat(\scr E)$.
	\begin{enumerate}
	\item The \emph{cyclic homology} of $\scr C$ is
	\[\HC(\scr C)=\abs{\scr C^\natural}_{h\T}\in\scr E.\]
	\item The \emph{negative cyclic homology} of $\scr C$ is
	\[\HN(\scr C)=\abs{\scr C^\natural}^{h\T}\in\scr E.\]
	\item If $\scr E$ is stable, the \emph{periodic cyclic homology} of $\scr C$ is
	\[\HP(\scr C)=\abs{\scr C^\natural}^{t\T}\in\scr E.\]
	\end{enumerate}
\end{definition}

Note that $\HC(\scr C)$ is simply the colimit of $\scr C^{\natural}\colon\Lambda^\op\to\scr E$. Note also that $\HC(\scr C)$, $\HN(\scr C)$, and $\HP(\scr C)$ depend only on the $\scr E$-enriched $\infty$-category presented by $\scr C$, since this is the case for $\HH(\scr C)$.

\begin{remark}
	There are several interesting refinements of the above definitions.	
	By definition, the invariants $\HC(\scr C)$, $\HN(\scr C)$, and $\HP(\scr C)$ depend only on the circle action on $\HH(\scr C)$. The \emph{topological cyclic homology} of $\scr C$ is a refinement of negative cyclic homology, defined when $\scr E$ is the $\infty$-category of spectra, which uses some finer structure on $\HH(\scr C)$.
	In another direction, additional structure on $\scr C$ can lead to $\HH(\scr C)$ being acted on by more complicated $\infty$-groups. For example, if $\scr C$ has a duality $\dagger$, then $\scr C^\natural$ extends to the dihedral category whose classifying space is $B\mathrm O(2)$. The coinvariants $\abs{\scr C^\natural}_{h\mathrm O(2)}$ are called the \emph{dihedral homology} of $(\scr C,\dagger)$.
\end{remark}

The previous definitions apply in particular when $\scr C$ has a unique object, in which case we may identify it with an $A_\infty$-algebra in $\scr E$. If $A$ is an $E_\infty$-algebra in $\scr E$, there is a more direct description of $\abs{A^\natural}$. In this case, $A^\natural$ is the underlying cyclic object of the cyclic $E_\infty$-algebra $\Lambda^0\tens A\in\PSh(\Lambda,\CAlg(\scr E))$, where $\Lambda^0$ is the cyclic set represented by $[0]\in\Lambda$ and $\tens$ is the canonical action of the $\infty$-category $\scr S$ of spaces on the presentable $\infty$-category $\CAlg(\scr E)$.
For any cyclic space $K\in\PSh(\Lambda)$, we clearly have $\abs{K\tens A}\simeq \abs{K}\tens A$. 
It follows that $\abs{A^\natural}\in\PSh(B\T,\scr E)$ is the underlying object of the $E_\infty$-algebra
\[
\abs{\Lambda^0}\tens A\simeq \T\tens A\in\PSh(B\T,\CAlg(\scr E)).
\]
In particular, $\HH(A)$, $\HN(A)$, and $\HP(A)$ inherit $E_\infty$-algebra structures from $A$. Their geometric interpretation is the following: if $X=\Spec A$, then $\Spec\HH(A)$ is the \emph{free loop space} of $X$ and $\Spec\HN(A)$ is the \emph{space of circles} in $X$. The cyclic homology $\HC(A)$ is a quasi-coherent sheaf on the space of circles in $X$, whose support is $\Spec\HP(A)$.

\section{Comparison with the classical definitions}
\label{sec:concrete}

Let $k$ be a discrete commutative ring and let $A$ be an $A_\infty$-algebra over $k$.
The cyclic and negative cyclic homology of $A$ over $k$ are classically defined via explicit bicomplexes. 
Let us start by recalling these definitions, following \cite[\S5.1]{Loday}. 

Let $M_\bullet$ be a cyclic object in an additive category $\scr A$.
The usual presentation of $\Lambda$ provides the face and degeneracy operators $d_i\colon M_n\to M_{n-1}$ and $s_i\colon M_n\to M_{n+1}$ ($0\leq i\leq n$), as well as the cyclic operator $c\colon M_n\to M_n$ of order $n+1$.
We define the additional operators
\begin{align*}
	b\colon M_n\to M_{n-1}, &\quad b=\sum_{i=0}^n(-1)^id_i,\\
	s_{-1}\colon M_n\to M_{n+1}, &\quad s_{-1}=cs_n,\\
	t\colon M_n\to M_n, &\quad t=(-1)^nc,\\
	N\colon M_n\to M_n,&\quad N=\sum_{i=0}^nt^i,\\
	B\colon M_n\to M_{n+1},&\quad B=(\id-t)s_{-1}N.
\end{align*}
We easily verify that $b^2=0$, $B^2=0$, and $bB+Bb=0$. In particular, $(M,b)$ is a chain complex in $\scr A$.
We now take $\scr A$ to be the category $\Ch_k$ of chain complexes of $k$-modules. Then $(M,b)$ is a (commuting) bicomplex and we denote by $(C_*(M),b)$ the total chain complex with
\[
C_n(M)=\bigoplus_{p+q=n}M_{p,q}, \quad b = b+(-1)^* d.
\]
We then form the (anticommuting) \emph{periodic cyclic bicomplex} $\BP(M)$:
\begin{tikzmath}
	\def\colsep{1.5em}
	\def\rowsep{1.5em}
	\diagram{& \vdots & \vdots & \vdots & \\
	\dotsb & C_2(M) & C_1(M) & C_0(M) & \dotsb \\
	\dotsb & C_1(M) & C_0(M) & C_{-1}(M) & \dotsb \\
	\dotsb & C_0(M) & C_{-1}(M) & C_{-2}(M) & \dotsb \\
	& \vdots & \vdots & \vdots & \\};
	\arrows
	(12) edge (22) (13) edge (23) (14) edge (24)
	(42) edge (52) (43) edge (53) (44) edge (54)
	(22) edge node[left]{$b$} (32) (23) edge node[left]{$b$} (33) (24) edge node[left]{$b$} (34)
	(32) edge node[left]{$b$} (42) (33) edge node[left]{$b$} (43) (34) edge node[left]{$b$} (44)
	(21-) edge[<-] (-22) (22-) edge[<-] node[above]{$B$} (-23) (23-) edge[<-] node[above]{$B$} (-24)
	(31-) edge[<-] (-32) (32-) edge[<-] node[above]{$B$} (-33) (33-) edge[<-] node[above]{$B$} (-34)
	(41-) edge[<-] (-42) (42-) edge[<-] node[above]{$B$} (-43) (43-) edge[<-] node[above]{$B$} (-44)
	(24-) edge[<-] (-25)
	(34-) edge[<-] (-35)
	(44-) edge[<-] (-45)
	;
\end{tikzmath}
with $\BP(M)_{p,q}=C_{q-p}(M)$.
Removing all the negatively graded columns, we obtain the \emph{cyclic bicomplex} $\BC(M)$; removing all the positively graded columns, we obtain the \emph{negative cyclic bicomplex} $\BN(M)$. Finally, we form the total complexes
\[
\Tot\BC,\, \Tot\BN,\, \Tot\BP\colon \PSh(\Lambda,\Ch_k)\to \Ch_k,
\]
where
\[
\Tot(B)_n=\colim_{r\to\infty}\prod_{p\leq r}B_{p,n-p}.
\]
These functors clearly preserve quasi-isomorphisms and hence induce functors
\[
\CC,\,\CN,\, \CP\colon \PSh(\Lambda,\Mod_k)\to\Mod_k,
\]
where $\Mod_k$ is the stable $\infty$-category of $k$-modules. There is a cofiber sequence
\[
\CC[1] \xrightarrow{B} \CN \to \CP,
\]
where the map ``$B$'' is induced by the degree $(0,1)$ map of bicomplexes $\BC(M) \to \BN(M)$ whose nonzero components are $B\colon C_{i-1}(M) \to C_i(M)$.

\begin{theorem}\label{thm:comparison}
	Let $k$ be a discrete commutative ring and $M\in\PSh(\Lambda,\Mod_k)$ a cyclic $k$-module. Then there are natural equivalences
	\[
	\abs{M}_{h\T}\simeq\CC(M)\text{,}\quad \abs{M}^{h\T}\simeq \CN(M)\text{, and}\quad \abs{M}^{t\T}\simeq \CP(M).
	\]
	In particular, if $\scr C$ is a $k$-linear $\infty$-category, then
	\[
	\HC(\scr C)\simeq\CC(\scr C^\natural)\text{,}\quad \HN(\scr C)\simeq \CN(\scr C^\natural)\text{, and}\quad \HP(\scr C)\simeq \CP(\scr C^\natural).
	\]
\end{theorem}

We first rephrase the classical definitions in terms of \emph{mixed complexes}, following Kassel \cite{Kassel}.
We let $k[\epsilon]$ be the differential graded $k$-algebra
\[\dotsb \to 0 \to k\epsilon \stackrel 0\to k\to 0 \to \dotsb,\]
which is nonzero in degrees $1$ and $0$.
The $\infty$-category $\Mod_{k[\epsilon]}$ is the localization of the category of differential graded $k[\epsilon]$-modules, also called mixed complexes, at the quasi-isomorphisms.
The functors 
\[k\tens_{k[\epsilon]}(\ph),\;\Hom_{k[\epsilon]}(k,\ph)\colon \Mod_{k[\epsilon]}\to\Mod_k\]
are related by a \emph{norm map}
\[
\nu_\epsilon\colon k[1]\tens_{k[\epsilon]}(\ph) \to \Hom_{k[\epsilon]}(k,\ph),
\]
induced by the $k[\epsilon]$-linear map $\epsilon\colon k[1]\to k[\epsilon]$.

We denote by
\[K\colon \PSh(\Lambda,\Mod_k)\to\Mod_{k[\epsilon]}\]
the functor induced by sending a cyclic chain complex $M$ to the mixed complex $(C_*(M),b,B)$.

\begin{lemma}\label{lem:Kassel}
	Let $M\in\PSh(\Lambda,\Mod_k)$. Then
	\[\CC(M)\simeq k\tens_{k[\epsilon]}K(M)\quad\text{and}\quad \CN(M)\simeq \Hom_{k[\epsilon]}(k,K(M)),\]
	and the norm map $\nu_\epsilon\colon k[1]\tens_{k[\epsilon]}K(M)\to \Hom_{k[\epsilon]}(k,K(M))$ is identified with $B\colon \CC(M)[1]\to \CN(M)$.
\end{lemma}

\begin{proof}
	We work at the level of complexes. Let $Qk$ be the nonnegatively graded mixed complex
	\[
	\dotsb\leftrightarrows k\epsilon\underset{0}{\overset{\epsilon}\leftrightarrows} k\underset{\epsilon}{\overset{0}\leftrightarrows} k\epsilon \underset{0}{\overset{\epsilon}\leftrightarrows} k.
	\]
	There is an obvious morphism $Qk\to k$ which is a cofibrant resolution of $k$ for the projective model structure on mixed complexes. By inspection, we have isomorphisms of chain complexes
	\[
	\Tot\BC(M)\simeq Qk\tens_{k[\epsilon]}(C_*(M),b,B)\quad\text{and}\quad \Tot\BN(M)\simeq \Hom_{k[\epsilon]}(Qk,(C_*(M),b,B)).
	\]
	This proves the first claim. 
	For $C$ a mixed complex, the norm map $\nu_\epsilon$ is modeled by the composition
	\[
	Qk[1]\tens_{k[\epsilon]} C \onto C/\mathrm{im}(\epsilon)[1] \xrightarrow \epsilon \ker(\epsilon|C) \into \Hom_{k[\epsilon]}(Qk,C).
	\]
	The last claim is then obvious.
\end{proof}

Let $k[\T]$ be the $A_\infty$-ring $k\tens\Sigma^\infty_+ \T$.
There is an obvious equivalence of $\infty$-categories
\[
\PSh(B\T,\Mod_k)\simeq \Mod_{k[\T]}
\]
that makes the following squares commute:
\begin{equation*}\label{eqn:square1}
\begin{tikzpicture}
	\diagram{\PSh(B\T,\Mod_k) & \Mod_{k} \\ \Mod_{k[\T]} & \Mod_{k}\rlap, \\};
	\arrows (21-) edge node[above]{$k\tens_{k[\T]}\ph$} (-22)
	(11-) edge node[above]{$(\ph)_{h\T}$} (-12) (11) edge node[left]{$\simeq$} (21) (12) edge[-,vshift=1pt] (22) edge[-,vshift=-1pt] (22);
\end{tikzpicture}
\qquad
\begin{tikzpicture}
	\diagram{\PSh(B\T,\Mod_k) & \Mod_{k} \\ \Mod_{k[\T]} & \Mod_{k}\rlap. \\};
	\arrows (21-) edge node[above]{$\Hom_{k[\T]}(k,\ph)$} (-22)
	(11-) edge node[above]{$(\ph)^{h\T}$} (-12) (11) edge node[left]{$\simeq$} (21) (12) edge[-,vshift=1pt] (22) edge[-,vshift=-1pt] (22);
\end{tikzpicture}
\end{equation*}

Since $\T$ is equivalent to the circle, $H_1(\T,\Z)$ is an infinite cyclic group.
Let $\gamma\in H_1(\T,\Z)$ be a generator. Sending $\epsilon$ to $\gamma$ defines an equivalence of augmented $A_\infty$-$k$-algebras $\gamma\colon k[\epsilon]\simeq k[\T]$, whence an equivalence of $\infty$-categories
\[
\gamma^*\colon \Mod_{k[\T]}\simeq \Mod_{k[\epsilon]}.
\]
Moreover, $\gamma^*$ identifies the norm maps $\nu_{\T}\colon \Sigma^{\mathfrak t}(\ph)_{h\T}\to (\ph)^{h\T}$ and $\nu_\epsilon \colon k[1]\tens_{k[\epsilon]}(\ph) \to \Hom_{k[\epsilon]}(k,\ph)$, provided that $\Sigma^{\mathfrak t}$ is identified with $\Sigma$ using the orientation of $\T$ given by $\gamma$ (since then the Poincaré duality isomorphism $k\simeq H_1(\T,k)$ sends $1$ to $\gamma$).

The main result of this note is that $K(M)$ is a model for $\abs{M}$. More precisely:

\begin{theorem}\label{thm:main}
	There exists a generator $\gamma\in H_1(\T,\Z)$ such that the following triangle commutes:
	\begin{tikzmath}
		\diagram{\PSh(\Lambda,\Mod_k) & \PSh(B\T,\Mod_k) \\ & \Mod_{k[\epsilon]}\rlap. \\};
		\arrows (11-) edge node[above]{$\abs{\ph}$} (-12) (12) edge node[right]{$\gamma^*$} node[left]{$\simeq$} (22) (11) edge node[below left]{$K$} (22);
	\end{tikzmath}
\end{theorem}

Theorem~\ref{thm:comparison} follows from Theorem~\ref{thm:main} and Lemma~\ref{lem:Kassel}. 
As an immediate corollary, we recover the following result of Dwyer and Kan \cite[Remark 6.7]{DK}:

\begin{corollary}
	The functor $K\colon \PSh(\Lambda,\Mod_k)\to \Mod_{k[\epsilon]}$ induces an equivalence of $\infty$-categories
	\[
	\PSh_\simeq(\Lambda,\Mod_k)\simeq \Mod_{k[\epsilon]}.
	\]
\end{corollary}

To prove Theorem~\ref{thm:main}, we consider the ``universal case'', namely the cocyclic cyclic $k$-module $k[\Lambda^\bullet]$. We have a natural equivalence
\[
M\simeq k[\Lambda^\bullet]\tens_\Lambda M,
\]
where 
\[
\tens_\Lambda\colon \Fun(\Lambda,\Mod_k)\times\PSh(\Lambda,\Mod_k)\to\Mod_k
\]
is the coend pairing. Similarly, we have
\[
\abs{M}\simeq \abs{k[\Lambda^\bullet]}\tens_\Lambda M\quad\text{and}\quad K(M)\simeq K(k[\Lambda^\bullet])\tens_\Lambda M,
\]
since both $\abs{\ph}$ and $K$ commute with tensoring with constant $k$-modules and with colimits (for $K$, note that colimits in $\Mod_{k[\epsilon]}$ are detected by the forgetful functor to $\Mod_k$).
Thus, it will suffice to produce an equivalence of cocyclic $k[\epsilon]$-modules
\begin{equation}\label{eqn:main}
	\gamma^*\abs{k[\Lambda^\bullet]}\simeq K(k[\Lambda^\bullet]).
\end{equation}

Let $k[u]$ denote the $A_\infty$-$k$-coalgebra $k\tens_{k[\epsilon]}k$. Note that a $k[u]$-comodule structure on $M\in\Mod_k$ is the same thing as map $M\to M[2]$.
The functor $k\tens_{k[\epsilon]}\ph\colon \Mod_{k[\epsilon]}\to \Mod_k$ factors through a fully faithful functor from ${k[\epsilon]}$-modules to $k[u]$-comodules:
\begin{tikzmath}
	\diagram{ & \Comod_{k[u]} \\ \Mod_{k[\epsilon]} & \Mod_k\rlap. \\};
	\arrows (21-) edge node[below]{$k\tens_{k[\epsilon]}\ph$} (-22) (12) edge node[right]{forget} (22) (21) edge[c->,dashed] (12);
\end{tikzmath}
To prove~\eqref{eqn:main}, it will therefore suffice to produce an equivalence of cocyclic $k[u]$-comodules
\begin{equation}\label{eqn:main2}
	k\tens_{k[\epsilon]}\gamma^*\abs{k[\Lambda^\bullet]}\simeq k\tens_{k[\epsilon]}K(k[\Lambda^\bullet]).
\end{equation}
Note that both cocyclic objects send all morphisms in $\Lambda$ to equivalences and hence can be viewed as functors $B\T\to\Comod_{k[u]}$. 

Let us first compute the left-hand side of~\eqref{eqn:main2}. The generator $\gamma$ induces an equivalence of coaugmented $A_\infty$-$k$-coalgebras $\check\gamma\colon k[u]\simeq k[B\T]$, whence an equivalence of $\infty$-categories
\[
\check\gamma^*\colon \Comod_{k[B\T]}\simeq \Comod_{k[u]}.
\]
We clearly have
\[
k\tens_{k[\epsilon]}\gamma^*\abs{k[\Lambda^\bullet]}\simeq \check\gamma^*\abs{k[\Lambda^\bullet]}_{h\T}.
\]
Now, $\abs{k[\Lambda^\bullet]}_{h\T}\simeq k[\abs{\Lambda^\bullet}_{h\T}]$, where $\abs{\Lambda^\bullet}_{h\T}$ is a $B\T$-comodule in $\Fun_{\simeq}(\Lambda,\scr S)\simeq \scr S_{/B\T}$. If $\pi^*\colon\scr S\to \scr S_{/B\T}$ is the functor $\pi^*X=X\times B\T$, then a $B\T$-comodule structure on $\pi^*X$ is simply a map $\pi^*X\to \pi^*B\T$, \ie, a map $X\times B\T\to B\T$ in $\scr S$. 
Here, $\abs{\Lambda^\bullet}_{h\T}$ is $\pi^*(*)\in\scr S_{/B\T}$ and its $B\T$-comodule structure $\sigma\colon\pi^*(*)\to\pi^*(B\T)$ is given by the identity $B\T\to B\T$. Applying $\check\gamma^*k[\ph]$, we deduce that the left-hand side of~\eqref{eqn:main2} is the constant cocyclic $k$-module $\underline{k}$ with $k[u]$-comodule structure given by the composition
\begin{equation}\label{eqn:comod1}
\underline k\stackrel\sigma\to\underline k[B\T]\stackrel{\check\gamma}\simeq \underline k[u].
\end{equation}
Note that equivalence classes of $k[u]$-comodule structures on $\underline{k}$ are in bijection with
\[
[\underline{k},\underline{k}[2]]\simeq H^2(B\T,k).
\]
Under this classification, \eqref{eqn:comod1} comes from an integral cohomology class, namely the image of the identity $B\T\to B\T$ under the isomorphism
\[
[B\T,B\T]\stackrel{\check\gamma}\simeq H^2(B\T,\Z).
\] 
In particular, it comes from a generator of the infinite cyclic group $H^2(B\T,\Z)$, determined by $\gamma$.
We must therefore show that the right-hand side of~\eqref{eqn:main2} is also equivalent to the constant cocyclic $k$-module $\underline k$ with $k[u]$-comodule structure classified by a generator of $H^2(B\T,\Z)$.

Recall that $K(k[\Lambda^\bullet])$ is the following mixed complex of cocyclic $k$-modules:
\[
\dotsb\rightleftarrows k[\Lambda_2]\underset{B}{\overset{b}{\rightleftarrows}}k[\Lambda_1]\underset{B}{\overset{b}{\rightleftarrows}} k[\Lambda_0].
\]
Consider the mixed complex $Qk$ from the proof of Lemma~\ref{lem:Kassel}, which can be used to compute $k\tens_{k[\epsilon]}\ph$ at the level of complexes. It comes with an obvious self-map $Qk\to Qk[2]$ which induces the $k[u]$-comodule structure on $k\tens_{k[\epsilon]}M$ for every mixed complex $M$. Let us write down explicitly the resulting chain complex $Qk\tens_{k[\epsilon]}K(k[\Lambda^\bullet])$ of cocyclic $k[u]$-comodules. It is the total complex of the first-quadrant bicomplex
\begin{tikzequation}\label{eqn:resolution}
	\def\colsep{1.5em}
	\def\rowsep{1.5em}
	\diagram{\vdots & \vdots & \vdots \\
	k[\Lambda_2] & k[\Lambda_1] & k[\Lambda_0] \\
	k[\Lambda_1] & k[\Lambda_0] & \\
	k[\Lambda_0]\rlap, & & \\};
	\arrows
	(11) edge (21) (12) edge (22) (13) edge (23)
	(21) edge node[left]{$b$} (31) (22) edge node[left]{$b$} (32)
	(31) edge node[left]{$b$} (41)
	(21-) edge[<-] node[above]{$B$} (-22) (22-) edge[<-] node[above]{$B$} (-23)
	(31-) edge[<-] node[above]{$B$} (-32)
	;
\end{tikzequation}
with $k[u]$-comodule structure induced by the obvious degree $(-1,-1)$ endomorphism $\delta$. 

\begin{proposition}\label{prop:resolution}
	The bicomplex~\eqref{eqn:resolution} is a resolution of the constant cocyclic $k$-module $\underline k$. Moreover, the endomorphism $\delta$ represents a generator of the invertible $k$-module $[\underline k,\underline k[2]]\simeq H^2(B\T,k)$.
\end{proposition}

\begin{proof}
	Let $K_{**}$ be the bicomplex~\eqref{eqn:resolution}, with the obvious augmentation $K_{**}\to\underline k$.
	For $M$ a cyclic object in an additive category, we define the operator $b'\colon M_n\to M_{n-1}$ by
	\[b'=b-(-1)^nd_n=\sum_{i=0}^{n-1}(-1)^id_i.\]
	 Let $L_{**}$ be the $(2,0)$-periodic first-quadrant bicomplex
	\begin{tikzmath}
		\def\rowsep{1.5em}
		\diagram{\vdots & \vdots & \vdots & \\
		k[\Lambda_2] & k[\Lambda_2] & k[\Lambda_2] & \dotsb \\
		k[\Lambda_1] & k[\Lambda_1] & k[\Lambda_1] & \dotsb \\
		k[\Lambda_0] & k[\Lambda_0] & k[\Lambda_0] & \dotsb \\
		};
		\arrows
		(11) edge (21) (12) edge (22) (13) edge (23)
		(21) edge node[left]{$b$} (31) (22) edge node[left]{$-b'$} (32) (23) edge node[left]{$b$} (33)
		(31) edge node[left]{$b$} (41) (32) edge node[left]{$-b'$} (42) (33) edge node[left]{$b$} (43)
		(21-) edge[<-] node[above]{$\id-t$} (-22) (22-) edge[<-] node[above]{$N$} (-23) (23-) edge[<-] (-24)
		(31-) edge[<-] node[above]{$\id-t$} (-32) (32-) edge[<-] node[above]{$N$} (-33) (33-) edge[<-] (-34)
		(41-) edge[<-] node[above]{$\id-t$} (-42) (42-) edge[<-] node[above]{$N$} (-43) (43-) edge[<-] (-44)
		;
	\end{tikzmath}
	with the obvious augmentation $L_{**}\to\underline k$,
	and let $M_{**}$ be the bicomplex obtained from $L_{**}$ by annihilating the even-numbered columns.
	Let $\phi\colon \Tot K_{**}\to \Tot L_{**}$ be the map induced by $(\id, s_{-1}N)\colon k[\Lambda_n]\to k[\Lambda_n]\oplus k[\Lambda_{n+1}]$, and let $\psi\colon \Tot L_{**}\to \Tot M_{**}$ be the map induced by $-s_{-1}N+\id\colon k[\Lambda_n]\oplus k[\Lambda_{n+1}]\to k[\Lambda_{n+1}]$. A straightforward computation shows that $\phi$ and $\psi$ are chain maps and that we have a commutative diagram with exact rows
	\begin{tikzequation}\label{eqn:exact}
		\def\colsep{1em}
		\diagram{0 & \Tot K_{**} & \Tot L_{**} & \Tot M_{**} & 0 \\
		0 & \underline k & \underline k & 0 & 0\rlap. \\};
		\arrows (11-) edge (-12) (12-) edge node[above]{$\phi$} (-13) (13-) edge node[above]{$\psi$} (-14) (14-) edge (-15)
		(21-) edge (-22) (22-) edge node[below]{$\id$} (-23) (23-) edge (-24) (24-) edge (-25)
		(12) edge (22) (13) edge (23) (14) edge (24);
	\end{tikzequation}
	From the identity $s_{-1}b'+b's_{-1}=\id$, we deduce that each column of $M_{**}$ has zero homology, and hence that that $\Tot M_{**}\simeq 0$. Next we show that each row of $L_{**}$ has zero positive homology, so that the homology of $\Tot L_{**}$ can be computed as the homology of the zeroth column of horizontal homology of $L_{**}$. This can be proved pointwise, so consider a part of the $n$th row evaluated at $[m]$:
\begin{equation}\label{eqn:cpx1}
\dotsb \to k[\Lambda(n,m)] \xrightarrow{\id-t} k[\Lambda(n,m)] \xrightarrow{N} k[\Lambda(n,m)] \to \dotsb.
\end{equation}
By the structure theorem for $\Lambda$, we have $\Lambda(n,m)=C_{n+1}\times\Delta(n,m)$, where $C_{n+1}$ is the set of automorphisms of $[n]$ in $\Lambda$. Thus, \eqref{eqn:cpx1} is obtained from the complex
\begin{equation}\label{eqn:cpx2}
	\dotsb \to k[C_{n+1}] \xrightarrow{\id-t} k[C_{n+1}] \xrightarrow{N} k[C_{n+1}] \to \dotsb.
\end{equation}
by tensoring with the free $k$-module $k[\Delta(n,m)]$, and we need only prove that~(\ref{eqn:cpx2}) is exact. Let \[x=\sum_{i=0}^nx_ic^{i}\in k[C_{n+1}].\]
Suppose first that $x(\id-t)=0$; then $x_i=(-1)^{ni}x_0$ and hence $x=x_0N$. Suppose next that $xN=0$, \ie, that $\sum_{i=0}^{n}(-1)^{ni}x_{n-i}=0$; putting $y_0=x_0$ and $y_i=x_i+(-1)^ny_{i-1}$ for $i>0$, we find $x=y(\id-t)$.
This proves the exactness of~\eqref{eqn:cpx2}, and also that the image of $\id-t\colon k[C_{n+1}]\to k[C_{n+1}]$ is exactly the kernel of the surjective map $k[C_{n+1}]\to k$, $x\mapsto \sum_{i=0}^n(-1)^{ni}x_{n-i}$. This map identifies the $0$th homology of the $n$th row of $L_{**}$ evaluated at $[m]$ with $k[\Delta(n,m)]$. Moreover, the vertical map $k[\Delta(n,m)]\to k[\Delta({n-1},m)]$ induced by $-b$ is the usual differential associated with the simplicial $k$-module $k[\Delta^m]$. This proves that $\Tot L_{\ast\ast}\to \underline k$ is a resolution of $k$.
From~\eqref{eqn:exact} we deduce that $\Tot K_{**}\to \underline k$ is a quasi-isomorphism.

To prove the second statement, we contemplate the complex $\Hom(\Tot K_{**},\underline k)$: it is the total complex of the bicomplex
\begin{tikzmath}
	\def\colsep{2em}
	\def\rowsep{1.5em}
	\diagram{
	& & k \\
	& k & k \\
	k & k & k \\
	\vdots & \vdots & \vdots \\};
	\arrows
	(13) edge node[left]{$0$} (23)
	(22) edge node[left]{$0$} (32) (23) edge node[left]{$\id$} (33)
	(31) edge node[left]{$0$} (41) (32) edge node[left]{$\id$} (42) (33) edge node[left]{$0$} (43)
	(22-) edge[<-] (-23)
	(31-) edge[<-] (-32) (32-) edge[<-] (-33)
	;
\end{tikzmath}
with trivial horizontal differentials and alternating vertical differentials. We immediately check that
\[
\Tot K_{**}\stackrel \delta\to (\Tot K_{**})[2]\to \underline k[2]
\]
is a cocycle generating the second cohomology module.
\end{proof}

It follows from Proposition~\ref{prop:resolution} that the right-hand side of~\eqref{eqn:main2} is the constant cocyclic $k$-module $\underline k$ with $k[u]$-comodule structure classified by $\delta\colon\underline k\to\underline k[2]$. Comparing with~\eqref{eqn:comod1} and noting that $\delta$ is natural in $k$, we deduce that Theorem~\ref{thm:main} holds by choosing $\gamma\in H_1(\T,\Z)$ to be the generator corresponding to $\delta\in H^2(B\T,\Z)$.

\providecommand{\bysame}{\leavevmode\hbox to3em{\hrulefill}\thinspace}

\end{document}